\documentclass[8pt]{article}
            
\usepackage{amsmath}               
\usepackage{amsfonts}              
\usepackage{amsthm}                
\usepackage[utf8]{inputenc}
\usepackage{leftidx}
\usepackage{tikz-cd}
\usepackage{mathtools}

\newtheorem{tw}{Theorem}[section]
\newtheorem{lem}{Lemma}[section]
\newtheorem{defi}{Definition}[section]

\newtheorem{rem}{Remark}[section]

\title{The basic $dd^{\mathcal{J}}$-lemma}
\author{Paweł Raźny}

\begin{document}
\maketitle

\begin{abstract}
The purpose of this short paper is to develop further the theory of transverse generalized complex structures (first introduced in \cite{Wade}). We focus on proving some equivalent conditions to the basic $dd^{\mathcal{J}}$-lemma in the spirit of \cite{Cav}. We justify our approach by describing the transverse symplectic structure in this language and relating the basic $dd^{\mathcal{J}}$-lemma to the surjectivity of the Lefschetz map. We also present a non-trivial example of a foliation endowed with a transverse generalized complex structure.
\newline\newline
$\bold{Classification}$: 53C12, 53D18, 55N99, 55T05.
\newline
$\bold{Keywords}$: basic cohomology, foliations, generalized complex structures, transverse geometric structures.
\end{abstract}

\section{Foliations}
In this section we are going to give a brief review of some basic facts concerning foliations and transverse structures.
\begin{defi} A codimension q foliation $\mathcal{F}$ on a smooth n-manifold M is given by the following data:
\begin{itemize}
\item An open cover $\mathcal{V}:=\{V_i\}_{i\in I}$ of M.
\item A q-dimensional smooth manifold $T_0$.
\item For each $V_i\in\mathcal{V}$ there is a submersion $f_i: V_i\rightarrow T_0$ with connected fibers (these fibers are called plaques).
\item For any intersection $V_i\cap V_j\neq\emptyset$ there exists a local diffeomorphism $\gamma_{ij}$ of $T_0$ such that $f_j=\gamma_{ij}\circ f_i$
\end{itemize}
The last condition ensures that plaques glue nicely to form a partition of M consisting of submanifolds of M of codimension q. This partition is called a foliation $\mathcal{F}$ of M and the elements of this partition are called leaves of $\mathcal{F}$.
\end{defi}
We call $T=\coprod\limits_{V_i\in\mathcal{V}}f_i(V_i)$ the transverse manifold of $\mathcal{F}$. The local diffeomorphisms $\gamma_{ij}$ generate a pseudogroup $\mathcal{H}$ of transformations on T (called the holonomy pseudogroup). The space of leaves $M\slash\mathcal{F}$ of the foliation $\mathcal{F}$ can be identified with $T\slash\mathcal{H}$.
\begin{defi}
 A smooth form $\omega$ on M is called basic if for any vector field X tangent to the leaves of $\mathcal{F}$ the following equality holds:
\begin{equation*}
i_X\omega=i_Xd\omega=0
\end{equation*}
Basic 0-forms will be called basic functions henceforth.
\end{defi}
Basic forms are in one to one correspondence with $\mathcal{H}$-invariant smooth forms on T. It is clear that $d\omega$ is basic for any basic form $\omega$. Hence, the set of basic forms of $\mathcal{F}$ (denoted $\Omega^{\bullet}(M\slash\mathcal{F})$) is a subcomplex of the de Rham complex of M. We define the basic cohomology of $\mathcal{F}$ to be the cohomology of this subcomplex and denote it by $H^{\bullet}(M\slash\mathcal{F})$. A transverse structure to $\mathcal{F}$ is an $\mathcal{H}$-invariant structure on T. For example:
\begin{defi}
$\mathcal{F}$ is said to be transversely symplectic if T admits an $\mathcal{H}$-invariant closed 2-form $\omega$ of maximal rank. $\omega$ is then called a transverse symplectic form. As we noted earlier, $\omega$ corresponds to a closed basic form of rank q on M (also denoted $\omega$).
\end{defi}
\begin{defi}
$\mathcal{F}$ is said to be transversely holomorphic if T admits a complex structure that makes all the $\gamma_{ij}$ holomorphic. This is equivalent to the existence of an almost complex structure $J$ on the normal bundle $N\mathcal{F}:=TM\slash T\mathcal{F}$ (where $T\mathcal{F}$ is the bundle tangent to the leaves) satisfying:
\begin{itemize}
\item $L_XJ=0$ for any vector field X tangent to the leaves.
\item if $Y_1$ and $Y_2$ are sections of the normal bundle then:
\begin{equation*}
 N_J(Y_1,Y_2):=[JY_1,JY_2]-J[Y_1,JY_2]-J[JY_1,Y_2]+J^2[Y_1,Y_2]=0
 \end{equation*}
where $[$ , $]$ is the bracket induced on the sections of the normal bundle.
\end{itemize}
\end{defi}
We also can define a special class of vector bundles and vector fields:
\begin{defi}
A vector bundle is called foliated if its transition functions are basic. Equivalently, it is a vector bundle on the transverse manifold on which the holonomy pseudogroup acts fiberwise linearly (by abuse of language we shall call such a bundle on the transverse manifold a foliated bundle as well). 
\end{defi}
\begin{defi}
A vector field $X\in\Gamma(TM)$ is called foliated if for all $Y\in\Gamma(T\mathcal{F})$ we have $[X,Y]\in\Gamma(T\mathcal{F})$. A section of the normal bundle is called foliated if it has a foliated vector field representative (equivalently all its representatives are foliated).
\end{defi}

\section{Transverse generalized complex structure}
Let us recall that for any manifold M there is a natural nondegenerate pairing $($ , $)_M$ on $TM\oplus T^*M$ defined at each point $x\in M$ by:
\begin{equation*}
(X_1+\alpha_1,X_2+\alpha_2)_M|_x:=\frac{1}{2}(\alpha_1(X_2)+\alpha_2(X_1))
\end{equation*}
for $X_i\in T_xM$ and $\alpha_i\in T^*_xM$. Let $\bold{T}(T_xM\oplus T_x^*M)$ denote the tensor algebra of $T_xM\oplus T_x^*M$ and let $I$ be the two-sided ideal in $\bold{T}(T_xM\oplus T_x^*M)$ generated by the elements of the form $v\otimes v -(v,v)_M|_x1_{\bold{T}(T_xM\oplus T_x^*M)}$ for $v\in T_xM\oplus T_x^*M$. The Clifford algebra:
\begin{equation*}
Cl(T_xM\oplus T_x^*M, (\text{ },\text{ } )_M|_x):=\bold{T}(T_xM\oplus T_x^*M)\slash I
\end{equation*}
acts on differential forms at $x$ by:
\begin{equation*}
(X+\alpha)\bullet\beta = i_X\beta+\alpha\wedge\beta
\end{equation*}
for any form $\beta$, covector $\alpha$ and vector $X$. There is also a generalization of the Lie bracket, called the Courant bracket, on $\Gamma(TM\oplus T^*M)$ defined by the formula:
\begin{equation*}
[X_1+\alpha_1,X_2+\alpha_2]:=[X_1,X_2]+i_{X_1}d\alpha_2-i_{X_2}d\alpha_1+\frac{1}{2}d(\alpha_2(X_1)-\alpha_1(X_2))
\end{equation*}
\begin{defi}
A generalized almost complex structure on M is an almost complex structure $\mathcal{J}$ on $TM\oplus T^*M$ orthogonal with respect to the natural pairing $($ , $)_M$ (i.e. $(Y_1,Y_2)_M=(\mathcal{J}Y_1,\mathcal{J}Y_2)_M$ for $Y_i\in \Gamma(TM\oplus T^*M)$). A generalized complex structure is a generalized almost complex structure satisfying the condition:
\begin{equation*}
N_{\mathcal{J}}(Y_1,Y_2):=[\mathcal{J}Y_1,\mathcal{J}Y_2]-\mathcal{J}[Y_1,\mathcal{J}Y_2]-\mathcal{J}[\mathcal{J}Y_1,Y_2]+\mathcal{J}^2[Y_1,Y_2]=0
\end{equation*}
for all $Y_i\in\Gamma(TM\oplus T^*M)$.
\end{defi}
From now on M will denote an n-dimensional manifold endowed with an $(n-2q)$-dimensional foliation $\mathcal{F}$ with transverse manifold T.
\begin{defi}
A transverse generalized complex structure on $(M,\mathcal{F})$ is a generalized complex structure on T which is invariant under the action of the holonomy pseudogroup. A foliated generalized complex structure is a generalized almost complex structure on the normal bundle constant along the leaves (i.e $L_X\mathcal{J}=0$ for $X\in\Gamma(T\mathcal{F})$) and integrable with respect to the Courant bracket modulo $T\mathcal{F}$ (i.e. $N_{\mathcal{J}}=0$ on basic forms and foliated vector fields in a small neighbourhood around any point $x\in M$). 
\end{defi}
\begin{rem} Foliated and transverse generalized complex structures are in a one to one correspondence.
\begin{proof}
This was proven in \cite{Wade} for transverse generalized almost complex structures. Showing that this correspondence induces a correspondence between integrable structures is a matter of local computations.
\end{proof}
\end{rem}
\begin{rem}
Note that the codimension of $\mathcal{F}$ has to be even due to the orthogonality condition (cf.\cite{Gul}).
\end{rem}
We will now slightly reformulate one of the main results of \cite{Wade}:
\begin{tw}
A transverse generalized complex structures on $\mathcal{F}$ is uniquely determined by a foliated subbundle L of $(N\mathcal{F}\oplus N^*\mathcal{F})\otimes\mathbb{C}$ which is involutive with respect to the Courant bracket modulo $T\mathcal{F}$ (i.e. in a small neighbourhood $V_x$ around any point $x\in M$ we have $[X,Y]\in\Gamma(L|_{V_x})$ for any $X,Y\in\Gamma(L|_{V_x})$ which are constant along the leaves), maximal isotropic with respect to the induced natural pairing on $(N\mathcal{F}\oplus N^*\mathcal{F})\otimes\mathbb{C}$ and satisfies $L\cap\bar{L}=0$. A subbundle satisfying the above conditions also determines a unique generalized complex structure.
\begin{proof}
Given a transverse generalized complex structure $\mathcal{J}$ one takes the i-eigen bundle of $\mathcal{J}$ on the transversal T. Since this bundle is invariant under the action of the holonomy pseudogroup it defines a foliated bundle on M. By \cite{Cav} (section 2.2) this bundle satisfies the properties stated in the theorem. Given a bundle L satisfying the properties stated in the theorem  it defines an invariant subbundle of $(TT\otimes T^*T)$ with analogous properties. This defines a transverse complex structure (again using section 2.2 of \cite{Cav}).
\end{proof}
\end{tw}
The integer $k=2q-dim_{\mathbb{C}}({\pi_{N_x\mathcal{F}\otimes\mathbb{C}}(L_x)})$ is called the type of the transverse generalized complex structure at point $x\in M$ (where $\pi_{N_x\mathcal{F}\otimes\mathbb{C}}$ is the projection onto $N_x\mathcal{F}\otimes\mathbb{C}$). A point $x$ is called regular if there is a neighbourhood of $x$ consisting of points at which the transverse generalized complex structure has type equal to its type at $x$. Let $B$ be a basic closed 2-form on $(M,\mathcal{F})$. If $\mathcal{J}$ is a generalized transverse complex structure, then we can define another transverse generalized complex structure by:
\begin{equation*}
\mathcal{J}^B_x:=
\left[ \begin{array}{cc}
1 & 0 \\
-B_x & 1 
\end{array} \right]
\mathcal{J}_x
\left[ \begin{array}{cc}
1 & 0 \\
B_x & 1 
\end{array} \right]
\end{equation*}
for any point $x\in M$ (in the above we treat $B_x$ as a linear map between the normal and conormal space at $x$). We call $\mathcal{J}^B$ a $B$-field transform of $\mathcal{J}$. We will now state a foliated version of the generalized Darboux theorem from \cite{Gul}.
\begin{tw} Let $(M,\mathcal{F})$ be a manifold endowed with a transversely generalized complex foliation of dimension $p=n-2q$ and let $x$ be a regular point of type k with respect to this structure. Then $x$ has a neighbourhood restricted to which $\mathcal{J}$ is equivalent, via a diffeomorphism, to a B-field transform of the standard generalized complex structure on $\mathbb{R}^p\times\mathbb{C}^k\times\mathbb{R}^{2q-2k}$ (for some basic closed 2-form B).
\begin{proof}
We can take any connected component of the transverse manifold containing an image of $x$. By the manifold version of the generalized Darboux theorem the image of $x$ has a neighbourhood $V$ for which the restriction of the generalized complex structure of the transverse manifold is equivalent via diffeomorphism to a B-field transform of the standard structure on $\mathbb{C}^k\times\mathbb{R}^{2q-2k}$ (for some closed form B). This generalized complex structure can be naturally extended to a foliated structure on $\mathbb{R}^p\times V$, which is diffeomorphic to a neighbourhood of $x$ (after the aforementioned extension B becomes basic). This structure gives the same structure as the original one on $V$, which means they are equivalent (via diffeomorphism).
\end{proof}
\end{tw}

\section{Example}
We are going to devote this section to constructing a non-trivial example of a transversely generalized complex foliation. Some examples of interest were given in \cite{Wade}. Simple examples include transversely holomorphic and transversely symplectic foliations (and their B-field transforms). To construct our example we are going to recall a foliated nilmanifold presented in \cite{Wol}.
\newline\indent We start by taking the group of upper-triangular complex matrices N and let $\mathcal{H}_s$ be the subgroup of N cosisting of the matrices of the form:
\begin{equation*}
\left[ \begin{array}{ccc}
1 & x_1+i(y_1+sy'_1) & x_3+sx'_3 + i(y_3+sy'_3)\\
0 & 1 & x_2+iy_2\\
0 & 0 & 1
\end{array} \right]
\end{equation*}
for $s\notin\mathbb{Q}$ and $x_j,x'_j,y_j,y'_j\in\mathbb{Z}$. The group $\mathcal{H}_s$ can be also considered as the $\mathbb{Z}^9$ subgroup of $\mathbb{R}^9$ with the following group operation:
\begin{eqnarray*}
(x_1,...,x_9)\square (a_1,...,a_9):=&(a_1+x_1,...,a_5+x_5,a_6+x_6+a_1x_4-a_2x_5,
\\
&a_7+x_7-a_3x_5,a_8+x_8+a_1x_5+a_2x_4,a_9+x_9+a_3x_4)
\end{eqnarray*}
There is also a submersion $u:(\mathbb{R}^9,\square)\rightarrow N$ given by:
\begin{equation*}
u(x_1,...,x_9):=(x_1+i(x_2+sx_3),x_4+ix_5,x_6+sx_7+i(x_8+sx_9))
\end{equation*}
This submersion is the identity when restricted to $\mathcal{H}_s$. Furthermore, it defines an $\mathcal{H}_s$-invariant foliation on $\mathbb{R}^9$, which in turn gives a foliation $\mathcal{F}$ on $M:=\mathbb{R}^9\slash\mathcal{H}_s$. We can define a transverse structure on $(M,\mathcal{F})$ by definining an $\mathcal{H}_s$-invariant structure on N. Since the standard complex structure on N is $\mathcal{H}_s$-invariant and the form $dx_2\wedge dy_2$ (where $dz_j=dx_j+idy_j$) is $\mathcal{H}_s$-invariant, the generalized complex structure in the chosen basis $(\frac{\partial}{\partial{x_j}},\frac{\partial}{\partial{y_j}},dx_j,dy_j)$ of $TN\oplus T^*N$ (for $j=1,2,3$) defined as:
\begin{equation*}
\mathcal{J}=
\left[ \begin{array}{rrrrrrrrrrrr}
0 & 0 & 0 & 1 & 0 & 0 & 0 & 0 & 0 & 0 & 0 & 0\\
0 & 0 & 0 & 0 & 0 & 0 & 0 & 0 & 0 & 0 & -1 & 0\\
0 & 0 & 0 & 0 & 0 & 1 & 0 & 0 & 0 & 0 & 0 & 0\\

-1 & 0 & 0 & 0 & 0 & 0 & 0 & 0 & 0 & 0 & 0 & 0\\
0 & 0 & 0 & 0 & 0 & 0 & 0 & 1 & 0 & 0 & 0 & 0\\
0 & 0 & -1 & 0 & 0 & 0 & 0 & 0 & 0 & 0 & 0 & 0\\

0 & 0 & 0 & 0 & 0 & 0 & 0 & 0 & 0 & 1 & 0 & 0\\
0 & 0 & 0 & 0 & -1 & 0 & 0 & 0 & 0 & 0 & 0 & 0\\
0 & 0 & 0 & 0 & 0 & 0 & 0 & 0 & 0 & 0 & 0 & 1\\

0 & 0 & 0 & 0 & 0 & 0 & -1 & 0 & 0 & 0 & 0 & 0\\
0 & 1 & 0 & 0 & 0 & 0 & 0 & 0 & 0 & 0 & 0 & 0\\
0 & 0 & 0 & 0 & 0 & 0 & 0 & 0 & -1 & 0 & 0 & 0
\end{array} \right]
\end{equation*}
is also $\mathcal{H}_s$-invariant. Integrability of this structure is obvious since this is the standard generalized complex structure on $\mathbb{C}^2\times\mathbb{R}^2$.

\section{$dd^{\mathcal{J}}$-lemma}
The results of section 2.2. of \cite{Cav} ensure that, there is a complex line subbundle $\Phi$ of the exterior algebra bundle $\bigoplus\limits_{i\geq 0} \wedge^i T^*M$ uniquely determined by the generalized complex structure. More precisely, for any point $x\in T$ it is uniquely determined by the property:
\begin{equation*}
L_x=\{(X+\alpha)\in(TT\oplus T^*T)\otimes\mathbb{C}\text{ }|\text{ } (X+\alpha)\bullet\Phi_x=0\}
\end{equation*}
At each point $x\in T$ this subbundle induces a decomposition of forms by:
\begin{equation*}
U^{q-k}_x:=(\wedge^k\bar{L}_x)\bullet\Phi_x
\end{equation*}
Each bundle $U^j$ is a foliated subbundle of $\bigoplus\limits_{i\geq 0} \wedge^i T^*M$ since L and $\bar{L}$ are both foliated. This decomposition induces a decomposition of the operator $d$:
\begin{equation*}
d=\partial+\bar{\partial} \quad \partial:\Gamma(U^j)\rightarrow \Gamma(U^{j+1}) \quad \bar{\partial}:\Gamma(U^j)\rightarrow \Gamma(U^{j-1})
\end{equation*}
Since each of the bundles $U^j$ is foliated the projections $\pi_j:\Omega^{\bullet}(T,\mathbb{C})\rightarrow \Gamma(U^j)$ take forms invariant under the action of the holonomy pseudogroup to forms invariant under the action of the holonomy pseudogroup. Hence, the same is true for the operators $\partial$ and $\bar{\partial}$. This allows us to consider these operators on basic forms. Furthermore, we can define the operator $d^{\mathcal{J}}:=i(\bar{\partial}-\partial)$. From now on by $d,d^{\mathcal{J}},\partial,\bar{\partial}$ we will always mean these operators restricted to basic forms. We shall also denote by $\Gamma_b(U^k)$ the set of basic sections of $U^k$.
\begin{defi}
We say that a foliation satisfies the basic $dd^{\mathcal{J}}$-lemma if:
\begin{equation*}
Im(dd^{\mathcal{J}})=Im(d)\cap Ker(d^{\mathcal{J}})=Ker(d)\cap Im(d^{\mathcal{J}})
\end{equation*}
This is by simple calculation equivalent to the basic $\partial\bar{\partial}$-lemma i.e.
\begin{equation*}
Im(\partial\bar{\partial})=Im(\partial)\cap Ker(\bar{\partial})=Ker(\partial)\cap Im(\bar{\partial})
\end{equation*}
\end{defi}
Throughout this section we are going to describe some equivalent conditions to the basic $dd^{\mathcal{J}}$-lemma in the spirit of \cite{Cav}.
\begin{tw}
Let $\mathcal{F}$ ba a transversely generalized complex foliation on a manifold M. $\mathcal{F}$ satisfies the basic $dd^{\mathcal{J}}$-lemma iff the inclusion 
\begin{equation*}
i:(Ker(d^\mathcal{J}),d)\rightarrow (\Omega^{\bullet}(M\slash\mathcal{F},\mathbb{C}),d)
\end{equation*}
induces an isomorphism in cohomology.
\begin{proof} First, let us assume that $\mathcal{F}$  satisfies the basic $dd^{\mathcal{J}}$-lemma. We denote the homomorphism induced in cohomology by $i^*$. To prove that this map is injective we shall show that it's kernel is trivial. We take $\alpha\in Ker(d^{\mathcal{J}})\cap Im(d)$ (cohomology classes of such forms constitute the kernel of $i^*$). By the basic $dd^{\mathcal{J}}$-lemma there exists a basic form $\beta$ such that $\alpha=dd^{\mathcal{J}}\beta$. This means that $\alpha$ is an image of a $d^{\mathcal{J}}$-closed basic form. Hence, the cohomology class of such a form is 0 in $(Ker(d^\mathcal{J}),d)$. This proves injectivity.
\newline\indent To prove surjectivity we are going to prove that for each closed basic form $\alpha$ there is a $d^{\mathcal{J}}$-closed basic form representing the same cohomology class. We define the form $\beta=d^{\mathcal{J}}\alpha$. Obviously, this form is $d$-closed and $d^{\mathcal{J}}$-exact. Hence, by the basic $dd^{\mathcal{J}}$-lemma there exists a basic form $\gamma$ satisfying $\beta=dd^{\mathcal{J}}\gamma$. Our $d^{\mathcal{J}}$-closed representative is the form $\alpha-d\gamma$.
\newline\indent On the other hand if $i$ is a quasi-isomorphism we take a form $\alpha$ such that $dd^{\mathcal{J}}\alpha=0$ i.e. $d^{\mathcal{J}}\alpha\in Im(d^{\mathcal{J}})\cap Ker(d)$. This means that $d\alpha$ represents the trivial class in  $(\Omega^{\bullet}(M\slash\mathcal{F},\mathbb{C}),d)$ and hence the trivial class in $(Ker(d^{\mathcal{J}}),d)$ (due to $i^*$ being an isomorphism). In other words, there exists $\beta\in Ker(d^{\mathcal{J}})$ such that $d\alpha=d\beta$. Using the fact that $i^*$ is an isomorphism we know that $[\alpha-\beta]$ has a $d^{\mathcal{J}}$-closed representative $\gamma$ i.e. $\alpha-\beta=\gamma+d\omega$. After applying $d^{\mathcal{J}}$ to this equality we get $d^{\mathcal{J}}\alpha=d^{\mathcal{J}}d\omega$. This proves that $Im(d^{\mathcal{J}}d)=Im(d^{\mathcal{J}})\cap Ker(d)$. The subsequent lemma finishes the proof.
\end{proof}
\end{tw}
\begin{lem}
Let $(M,\mathcal{F})$ be a foliated manifold endowed with a transverse generalized complex structure $\mathcal{J}$. Then the following conditions are equivalent:
\begin{enumerate}
\item $Im(dd^{\mathcal{J}})=Im(d^{\mathcal{J}})\cap Ker(d)$
\item $Im(dd^{\mathcal{J}})=Im(d)\cap Ker(d^{\mathcal{J}})$
\end{enumerate}
\begin{proof} We shall prove only the implication "$1.\implies 2.$" as the proof of the converse is analogous. We can define an action of $\mathcal{J}$ on $U^k$ as multiplication by $i^k$ (see \cite{Cav} for motivation). These actions can be combined into an action of $\mathcal{J}$ on $\Omega^{\bullet}(M\slash\mathcal{F},\mathbb{C})$. Furthermore, this action takes invariant forms on the transverse manifold into invariant forms, because the bundles $U^k$ are complex and foliated. The equation $d^{\mathcal{J}}=\mathcal{J}^{-1}d\mathcal{J}$ holds on the transverse manifold (by e.g. \cite{Cav}) and hence also holds for basic forms. Using this formula it is easy to see that if $\alpha\in Im(d)\cap Ker(d^{\mathcal{J}})$  then $\mathcal{J}\alpha\in Im(d^{\mathcal{J}})\cap Ker(d)$. Due to $1.$ we have:
\begin{equation*}
\alpha=-\mathcal{J}dd^{\mathcal{J}}\beta=-\mathcal{J}d\mathcal{J}^{-1}d\mathcal{J}\beta=dd^{\mathcal{J}}(-\mathcal{J}\beta)
\end{equation*}
Hence, $\alpha\in Im(dd^{\mathcal{J}})$.
\end{proof}
\end{lem}
\begin{tw}
Let $\mathcal{F}$ be a transversely generalized complex foliation on a manifold M. $\mathcal{F}$ satisfies the basic $dd^{\mathcal{J}}$-lemma iff the projection:  
\begin{equation*}
p:(\Omega^{\bullet}(M\slash\mathcal{F},\mathbb{C}),d^{\mathcal{J}})\rightarrow (\Omega^{\bullet}(M\slash\mathcal{F},\mathbb{C})\slash Im(d),d^{\mathcal{J}})
\end{equation*}
induces an isomorphism in cohomology (here we consider $\Omega^{\bullet}(M\slash\mathcal{F})$ as a $\mathbb{Z}_2$-graded algebra with the gradation given by the parity of the degree of forms).
\begin{proof}
As in the previous theorem we denote by $p^*$ the homomorphism induced in cohomology. Since the kernel of $p^*$ is represented by forms in $Im(d)\cap Ker(d^{\mathcal{J}})$ it is easy to see that injectivity of $p^*$ is equivalent to $Im(d)\cap Ker(d^{\mathcal{J}})=Im(d)\cap Im(d^{\mathcal{J}})$. A form $\alpha$ represents a cohomology class of $(\Omega^{\bullet}(M\slash\mathcal{F},\mathbb{C})\slash Im(d),d^{\mathcal{J}})$ iff $d^{\mathcal{J}}\alpha\in Im(d)$. Surjectivity of $p^*$ is equivalent to the existence of a $d^{\mathcal{J}}$-closed form in every class of $(\Omega^{\bullet}(M\slash\mathcal{F},\mathbb{C})\slash Im(d),d^{\mathcal{J}})$. This in turn is equivalent to $Im(dd^{\mathcal{J}})=Im(d)\cap Im(d^{\mathcal{J}})$. This together with the previous lemma ends the proof.
\end{proof}
\end{tw}
We shall end this section by showing the correlation between the basic $dd^{\mathcal{J}}$-lemma and the existence of a decomposition in cohomology induced by the subbundles $U^k$. To this end we introduce the canonical spectral sequence as the spectral sequence associated to the double complex:
\begin{equation*}
(K^{p,q}:=\Gamma_b(U^{p-q}),\bar{\partial}:K^{\bullet,\bullet}\rightarrow K^{\bullet ,\bullet +1},\partial:K^{\bullet,\bullet}\rightarrow K^{\bullet +1,\bullet}) 
\end{equation*}
given by the filtration induced by the first degree $p$ (i.e such that the first page is the cohomology of the double complex with respect to $\partial$).
\begin{tw}
Let $\mathcal{F}$ be a transversely generalized complex foliation on a manifold M. The following conditions are equivalent:
\begin{enumerate}
\item $\mathcal{F}$ satisfies the basic $dd^{\mathcal{J}}$-lemma
\item The canonical spectral sequence degenerates at the first page and the subbundles $U^k$ induce a decomposition in cohomology.
\end{enumerate}
\begin{proof} We only need to prove that the basic $dd^{\mathcal{J}}$-lemma implies that the bundles $U^k$ induce a decomposition in cohomology. Other required implications are a consequence of Theorem 5.17 in \cite{Del} applied to our case. Any given cohomology class has a representative $\alpha$ which is $d^{\mathcal{J}}$-closed (by Theorem 4.1.). This means that $\alpha$ is both $\partial$-closed and $\bar{\partial}$-closed. Hence, all $\alpha_k:=\pi_k(\alpha)$ are $\partial$-closed and $\bar{\partial}$-closed and consequently $d$-closed. Furthermore, if the chosen class is zero, then the form $\alpha$ is exact and $d^{\mathcal{J}}$-closed. By the basic $dd^{\mathcal{J}}$-lemma $\alpha=dd^{\mathcal{J}}\beta$ for some form $\beta$. Since the degree of $dd^{\mathcal{J}}$ is zero with respect to the grading given by $U^k$ this implies that $\alpha_k=dd^{\mathcal{J}}\beta_k$. Hence, each $\alpha_k$ is exact.
\end{proof}
\end{tw}

\section{Transversely symplectic foliations}
In this section we are going to describe the transverse symplectic structure and the basic $dd^{\Lambda}$-lemma in the language of transverse generalized complex structures. Everything in this section except for the final theorem is a simple consequence of analogous statements on the transverse manifold. Given a transversely symplectic foliation $(M,\mathcal{F},\omega)$ we can define a transverse generalized complex structure on $\mathcal{F}$:
\begin{equation*}
\mathcal{J}:=
\left[ \begin{array}{cc}
0 & -\omega^{-1}\\
\omega & 0
\end{array} \right]
\end{equation*}
One can define a transverse symplectic star operator $*_s$ by defining it on the tangent manifold. With the help of the symplectic star we can define the following important operators:
\begin{equation*}
L(\alpha):=\omega\wedge\alpha \quad \Lambda(\alpha):=\star_s L\star_s (\alpha) \quad d^{\Lambda}\alpha:=d\Lambda(\alpha)-\Lambda d(\alpha)
\end{equation*}
for a basic form $\alpha$ (the operators $L$ and $\Lambda$ are also well defined pointwise). The operator $d^{\Lambda}$ coincides with $d^{\mathcal{J}}$ due to the fact that their analouges coincide on the transverse manifold. Using the introduced notation established we can describe $U^k$ with the help of the operators above:
\begin{equation*}
U^k_x=\{e^{i\omega} e^{\frac{\Lambda}{2i}} \alpha\text{ } |\text{ } \alpha\in\Omega^k_x(T,\mathbb{C})\}
\end{equation*}
We can also compute $\partial$ and $\bar{\partial}$:
\begin{equation*}
\partial e^{i\omega}e^{\frac{\Lambda}{2i}}\alpha=e^{i\omega}e^{\frac{\Lambda}{2i}}d\alpha
\quad
\bar{\partial}e^{i\omega}e^{\frac{\Lambda}{2i}}\alpha=-e^{i\omega}e^{\frac{\Lambda}{2i}} d^{\Lambda}\alpha
\end{equation*}
With this description we get the following theorem as a corollary from the discussion in the previous section:
\begin{tw}
With notation as above the following conditions are equivalent:
\begin{enumerate}
\item $\mathcal{F}$ satisfies the basic $dd^{\Lambda}$-lemma.
\item $i:(Ker(d^{\Lambda}),d)\rightarrow (\Omega^{\bullet}(M\slash\mathcal{F}),d)$ induces an isomorphism in cohomology.
\item The subbundles $U^k$ induce a decomposition in cohomology.
\item $p:(\Omega^{\bullet}(M\slash\mathcal{F},\mathbb{C}),d^{\Lambda})\rightarrow (\Omega^{\bullet}(M\slash\mathcal{F},\mathbb{C})\slash Im(d),d^{\Lambda})$ induces an isomorphism in cohomology.
\end{enumerate}
Furthermore, each of them implies that the Lefschetz map
\begin{equation*}
L^k:H^{q-k}(M\slash\mathcal{F})\rightarrow H^{q+k}(M\slash\mathcal{F})
\end{equation*}
is surjective.
\begin{proof} The first three conditions are equivalent due to the discussion in the previous section and the fact that for a symplectic manifold the canonical spectral sequence always degenerates at the first page as was proven in \cite{My}. Since, $d^{\Lambda}$ is of degree $-1$ we can repeat the proof of Theorem 4.2 with the $\mathbb{Z}$-grading to prove that $1.$ and $4.$ are equivalent. Condition $2.$ implies that every cohomology class has $d^{\Lambda}$-closed representative. This property is equivalent to the surjectivity of the Lefschetz map as was shown in \cite{BC}. 
\end{proof}
\end{tw}
Note that the surjectivity of the Lefschetz map does not imply the basic $dd^{\Lambda}$-lemma in general. A simple counterexample is $\mathbb{R}^2$ with the standard symplectic structure and the unique foliation of dimension $0$. A stronger version of this theorem can be proven for transversely symplectic, homologically orientable, Riemannian foliations on compact manifolds. Under these assumptions surjectivity of the Lefschetz map is equivalent to it being an isomorphism (by the basic version of Poincare duality) and this in turn is equivalent to the $dd^{\Lambda}$-lemma (cf. \cite{My}).

\end{document}